\documentclass[a4paper,10pt]{article}
\usepackage[cp1251]{inputenc}
\usepackage[english]{babel}

\usepackage{amssymb}\usepackage{amsmath}\usepackage{amsthm}
\usepackage{epsfig}
\def\be{\begin{equation}}
\def\ee{\end{equation}}
\newtheorem{thrm}{\bf Theorem}
\newtheorem{lmm}{\bf Lemma}
\newtheorem*{dfn}{\bf Definition}
\newtheorem*{rmk}{\bf Remark}
\newtheorem*{xmpl}{\bf Example}

\begin{document}

\title {How many Zolotar\"ev fractions are there?} 
\author{\copyright 2015 ~~~~A.B.Bogatyr\"ev
\thanks{Supported by RFBR grants 13-01-00115, 13-01-12417ofi-m2 and RAS Program 
"Modern problems of theoretical mathematics"}} 
\date{} 
\maketitle

E.I.Zolotar\"ev has found the best uniform rational approximation for the 
$Sign$ function on two real intervals separated by zero \cite{Z1, A1}. Known today as Zolotar\"ev fractions, 
those functions possess lots of interesting properties which in particular led to their applications in electrical engineering
(Cauer or elliptic filters). Despite rather involved representation (see below), the itself construction of Zolotar\"ev fractions
seems very natural and it is not a wonder that they were rediscovered several times \cite{Todd, Ng, Luto} since 1877. 

The parametric formulas for classical degree $n$ Zolotar\"ev fractions $Z_n(x|\tau)$  are given in terms of 
elliptic sine $sn(\cdot|\tau)$ and complete elliptic integral $K(\tau)$ of modulus $\tau\in i\mathbb{R}_+$  \cite{A1,A2,B2} 
\be
Z_n(x_{n\tau}(u)|\tau):=x_\tau(u),  \qquad u\in \mathbb{C}, 
\label{Zn}
\ee
where $x_\tau(u)=sn(K(\tau)u|\tau)$. At the first glance it is not clear that $Z_n$ is a rational function of 
$x$, same difficulty we encounter even for Chebyshev polynomials $T_n(cos(u))=cos(nu)$. 

The geometric interpretation (most likely not new) for the above formula (\ref{Zn}) leads to a natural enumerative problem for the arising construction. The solution to the latter problem is given in this note.

\section{Twisted Zolotar\"ev fractions and their classes}
One of the charachteristic features of the above functions $Z_n$ we take as a 

\begin{dfn}
A rational function $R(x)$ we call twisted Zolotar\"ev fraction iff all of its critical points (including those at infinity if any)
are simple with the values in a 4-element set. Two fractions are in the same class iff they differ by a pre- and post- 
compositions with complex linear fractional functions. 
\end{dfn}

\begin{rmk} 
1) The above definition of (twisted) Zolotar\"ev fractions is not standard and it is a replica of $(g=1)$ extremal rational functions \cite{B1}. Classical fractions (\ref{Zn}) surely meet this definition when $n>2$: their critical points are exactly inner alternation points corresponding to the mentioned above uniform approximation problem \cite{A1}.

2)The global change of dependent and independent variables of a rational function by elements of $PSL_2(\mathbb{C})$
cannot be considered as a significant modification of this function. In particular, it keeps intact all its discrete invariants:
the degree, the number of critical points and their indicies, the monodromy \cite{ZL}. We shall see below that recently introduced 
"Elliptic rational functions" \cite{Luto} and "Chebyshev-Blaschke products" \cite{Ng} fall  into the same class with classical Zolotar\"ev fractions.
\end{rmk}

The sets $B$ of four critical values of twisted Zolotar\"ev fractions from the same class make up an orbit of the 
$PSL_2(\mathbb{C})$-action on 4-tuples of points. Those orbits are distinguished by the value of the so called $j$-invariant \cite{A2,TM}:
\be
j:=\frac4{27}\frac{(\lambda^2-\lambda+1)^3}{\lambda^2(1-\lambda)^2},
\ee
where $\lambda$ is a cross-ratio of the four points containing in $B$ and the right hand side of the 
expression is the principal invariant of the so called unharmonic group. The latter interchanges the values of the cross ratio
$$
\lambda,~ 1-\lambda,~ \frac\lambda{\lambda-1},~ \frac1{1-\lambda},~ \frac1{\lambda},~ \frac{\lambda-1}\lambda,
$$
under the permutation of points in $B$.

{\bf The main purpose} of this note is to find the number of classes of twisted Zolotar\"ev fractions,
given the degree and the (finite) value of $j$-invariant.

\section{Construction of twisted Zolotar\"ev fractions}
\label{LatCnstr}
Consider  rank two lattice $ L$ in the complex plane of variable $u$; the group of translations of the plane by the elements of the lattice we designate by the same letter $ L$. Let $ L^+$ be the extension of the translation group by the  central symmetry: $u\to -u$. The extended group acts discontinuously in the complex plane, so its orbit space is well defined and carries natural complex structure
$$
{\mathbb C}/ L^+={\mathbb C}P^1.
$$
We can introduce a global coordinate on this Riemann sphere, say 
\be
\label{xL}
x_L(u)=\wp(u| L):=u^{-2}+\sum_{0\neq v\in L}((u-v)^{-2}-v^{-2}),
\ee
transition to any other global coordinate is given by a linear fractional map.

Once we have a full rank sublattice $M$ of $L$, the group $M^+$ is a subgroup of $ L^+$
and any orbit of $ M^+$ is contained in the orbit of $L^+$. Therefore we have a holomorphic mapping
from one sphere to the other:
\be
{\mathbb C}/ M^+\to {\mathbb C}/ L^+, 
\qquad u\in \mathbb{C},
\ee 
which becomes a rational function once we fix complex coordinate on each sphere. Thus we obtain a 
degree $|L:M|$ rational function $R_{L:M}(x)$:
\be
R_{L:M}(x_M(u)):= x_L(u),   
\qquad  x_M(u):= \wp(u| M),
\qquad u\in \mathbb{C}.
\label{R}
\ee

\begin{xmpl}
(1)To get modulus  $\tau\in i\mathbb{R}_+$ classical Zolotar\"ev fraction  we take $ L=L(\tau):=Span_\mathbb{Z}\{4, 2\tau\}$ and
$ M:=Span_\mathbb{Z}\{4, 2n\tau\}=:L(n\tau)$, then $R_{ L: M}(x)$ and $Z_n(x|\tau)$ lie in the same class, that is  
they differ by pre- and post- compositions with linear fractional functions. Indeed,  $x_\tau(u+1):=sn(K(\tau)(u+1)|\tau)$ 
is a degree two even elliptic function of $u$ with the period lattice $L(\tau)$, therefore it can be taken as a global coordinate in the orbit space of the group $L^+(\tau)$ instead of Weierstrass function (\ref{xL}). Classical Zolotar\"ev fractions are therefore described by the construction (\ref{R}).

(2) Elliptic rational functions $Y_n(x|\tau)$  also dubbed as Chebyshev-Blaschke products are defined as follows \cite{Luto, Ng}:
$$
Y_n(y_{1,\tau}(u)|\tau):=y_{n,\tau}(u), \qquad u\in\mathbb{C},
$$
where $ y_{n,\tau}(u):=\sqrt{k(n\tau)}cd(2K(n\tau)nu|n\tau)$. 
The latter function $ y_{n,\tau}$ is a degree two even elliptic function of $u$  with the period lattice $L_n(\tau):=Span_\mathbb{Z}\{2/n, \tau\}$ which may be taken as a global coordinate in $\mathbb{C}/L_n^+$. Hence, the function $Y_n(\cdot|\tau)$ fits the description (\ref{R}) with $L=L_n(\tau)$ and $M=L_1(\tau)$. 

(3) When $M=nL$, function (\ref{R}) coinsides with the so called flexible Latt\'es map \cite{JMil} of degree $n^2$.
In general, the construction does not reduce to a Latt\'es map.
\end{xmpl}

\begin{thrm}
\label{ZfrEqLat}
Construction (\ref{R}) gives 1-1 correspondence between
classes  of twisted Zolotar\"ev fractions and the pairs of
embedded rank two lattices $M\subset L\subset\mathbb{C}$ modulo multiplication by nonzero complex number, with two exceptions:
(1) $|L:M|=2$ and (2) $L=2M$.
\end{thrm}
{\bf Proof. (1).} 
Let $M\subset L$ be two lattices in the complex plane and $R_{L:M}$ be a corresponding rational function (\ref{R}).
Let us check that it has only simple critical points with the values in a 4-element set.
Indeed, critical points of $x_L(u)$ correspond to fixed points of the $L^+$ action on the plane. They are simple and make up the lattice $\frac12 L$.  Critical values of the projection $x_L(u)$ make up a 4-element set since $|\frac12 L:L|=4$. 
Critical points of the function $R_{L:M}$ correspond to the set
$\frac12 L\setminus \frac12 M$ in the complex plane.

Let us count how many different orbits of $L^+$ are there in the set $\frac12 L\setminus\frac12 M$.
One can show (use Schmitt normal form of integer matrices \cite{Newman}) that there are bases in the lattice and its sublattice related by an integer diagonal matrix. Without loss of generality we assume that $L=Span_\mathbb{Z}(1,\tau)$ and $M=Span_\mathbb{Z}(p,q\tau)$ with positive integers $p\ge q$ and $Im~\tau\neq0$. When $p\ge3$ no one of four points $u=\frac12$, $1$, $\frac12+\frac12\tau$, $1+\frac12\tau$ belongs to the lattice $M/2$, and no two of them are in the same orbit of $L^+$, therefore $R_{L:M}(x)$ has exactly four critical values.

Otherwise, if $p=q=2$, we have $M=2L$ and the set $(\frac12 L\setminus\frac12 M)/L^+$ has only three points represented e.g. by 
$u=\frac12$, $\frac12+\frac12\tau$, $\frac12\tau$. The corresponding function $R_{L:M}(x)$  falls into the class 
of one of the simplest Belyi functions: $~x^2+x^{-2}~$, the Hauptmodul for the Klein's quadratic group.
The remaining case $p=2$, $q=1$ corresponds to the index $2$ sublattice and  $R_{L:M}(x)$ has just two critical points as well as values, therefore it is in the class of the quadratic function $x^2$.

In case we rescale the lattices, the rational function $R_{L:M}$ changes in a predictable way since
$$
\wp(cu|cL)=c^{-2}\wp(u,L), \qquad c\in \mathbb{C}^\times.
$$
In particular, $R_{L:M}(x)$ remains within the same class of Zolotar\"ev fractions.

{\bf (2)}. Conversely, let $R(x)$ be a degree $n$ twisted Zolotar\"ev fraction with the set $B$ of four distinct critical values.
We reveal two nonexceptional lattices $M\subset L$ such that $R(x)$ is in the same class as $R_{L:M}(x)$.

Riemann-Hurwitz formula suggests that there are exactly four noncritical points in the pre-image  $R^{-1}B$.
We designate this set as $B_1$. Let $T$ (resp. $T_1$) be a torus doubly covering Riemann sphere
with branching at the points of $B$ (resp. $B_1$). In what follows, all the objects related to the torus $T_1$ we supply 
with the subindex $1$ by default. When the fraction $R$ changes in its class, the branching sets $B$ and $B_1$ 
are subjected to the action of $PSL_2(\mathbb{C})$, so  we assume wlog that $\infty\not\in B, ~B_1$. The 
affine part of the elliptic curve $T$ is now given by the equation: 
\be
(x,w)\in\mathbb{C}^2:
\qquad w^2=\prod\limits_{e\in B} (x-e),
\label{T}
\ee
and the curve admits an involution $J(x,w):=(x,-w)$ with four fixed points which project to the set $B$.  
We claim that there exists an isogeny $\tilde{R}: ~T_1\to T$
of the tori which commutes with the involutions $J,J_1$ and it's descent to the bases $T_1/J_1$ and $T/J$
coinsides with $R(x)$. The isogeny is a (essentially unique) rational map

\be
\tilde{R}: T_1\ni(x_1,w_1)~~ \to ~~(R(x_1), w_1S(x_1)):=(x,w)\in T,
\label{liftR}
\ee
where $S(\cdot)$ is a rational function with simple zeroes at the critical points of $R(\cdot)$, double zero at infinity and double poles at the poles of $R(\cdot)$.

Let us recall some facts about universal covering $\tilde{T}$ of the torus. The latter is the space of paths on $T$ starting at the marked point modulo their homotopy. To introduce a global complex coordinate on the universal covering we choose a (unique up to nonzero factor) holomorphic differential on the torus, for instance $d\eta=w^{-1}dx$  in the algebraic model (\ref{T}) of the torus. Now we assign a complex number $u[\gamma]:=\int_\gamma~ d\eta$ to a path $\gamma$. The torus itself is a factor of the complex plane by the group of cover transformations which is realized as translations by the elements of the 
lattice 
$$
L:=\int_Cd\eta,
\qquad C\in H_1(T,\mathbb{Z}).
$$
If the marked point on the torus is fixed by the involution, then the lift of $J$ to $\tilde{T}$ is just a central symmetry:
$$
u[J\gamma]=\int_{J\gamma}d\eta=\int_\gamma J^*d\eta=-\int_\gamma d\eta=-u[\gamma].
$$ 
The natural projection $\mathbb{C}\to\mathbb{C}P_1=T/J=\mathbb{C}/L^+$
is a superposition of a linear fractional map (which we do not specify here) and the function $x_L(u):=\wp(u|L)$.

The torus $T$ and its covering torus $T_1$ share the universal covering in a natural way: back and forth mappings 
$\tilde{T}\leftrightarrow\tilde{T}_1$ are given by the projection and lifting of paths,
provided the marked points on the tori agree. For the identified points of $\tilde{T}$ and $\tilde{T}_1$ to have the same complex coordinate, we take the pullback of the  holomorphic differential $d\eta$ as a distinguished differential $d\eta_1$ on $T_1$. Now the identified paths $\gamma_1\subset T_1$ and $\tilde{R}\gamma_1:=\gamma\subset T$ will share the coordinate: 
$$
u_1[\gamma_1]:=\int_{\gamma_1}d\eta_1=\int_{\gamma_1}\tilde{R}^*d\eta=\int_{\tilde{R}\gamma_1}d\eta:=u[\gamma].
$$
With this choice of the global coordinate on $\tilde{T}_1$, the covering torus will be a factor of the complex plane  by the sublattice of $L$:
$$
M=L_1:=
\int_{H_1(T_1,\mathbb{Z})}~d\eta_1=
\int_{H_1(T_1,\mathbb{Z})}~\tilde{R}^*d\eta=
\int_{\tilde{R}H_1(T_1,\mathbb{Z})}~d\eta
\subset\int_{H_1(T,\mathbb{Z})}~d\eta:=L.
$$

The consistent choice of the marked points on the tori: $(e_1,0)\in T_1$, $e_1\in B_1$ and 
respectively $(e,0)\in T$, $e:=R(e_1)\in B$ leads to a commutative square:

\be
\begin{array}{rcl} 
\tilde{T}_1=\mathbb{C}&\stackrel{id}{\longrightarrow}&\mathbb{C}=\tilde{T}\\

l_1\circ x_M(u)\downarrow~~~ && ~~~\downarrow l\circ x_L(u)\\

\mathbb{C}/M^+=\mathbb{C}P^1 & \stackrel{R(x)}{\longrightarrow} &\mathbb{C}P^1=\mathbb{C}/L^+\\

\end{array}
\ee

with linear fractional functions $l_1,~l$, wherefrom the representation (\ref{R}) follows:
$$
R=l\circ x_L\circ x_M^{-1}\circ l_1^{-1}.
$$

The arising pair of lattices $M\subset L$ are not exceptional since they generate a function $R_{L:M}(x)$ with 
four critical values, as many as $R(x)$ has. Another choice of the distinguished differential 
$d\eta$ on $T$ leads to rescaling/rotation of the lattices. This ends the proof.

\section{Zolotar\"ev and Chebyshev-Blaschke are classmates}
Classical Zolotar\"ev fractions and Chebyshev-Blaschke products as it follows from Examples (1) and (2) of Sect. \ref{LatCnstr} are determined by proportional pairs of embedded lattices $(L(\tau), L(n\tau))=n\tau~(L_n(\tau_1), L_1(\tau_1))$  with $\tau_1=-4/(n\tau)$. By Theorem \ref{ZfrEqLat},  two mentioned rational functions differ by normalizations of dependent and independent variables only. We show this in a strightforward way below.

First we establish the key identity 
\begin{equation}
cd(2K(\tau)u|\tau)=l_\tau\circ cd(iK'(\tau_0)u|\tau_0),
\label{ident}
\end{equation}
where $\tau_0=-4/\tau$, $K(\tau), K'(\tau)$ are complete elliptic integrals of modulus $\tau$ (they may be defined in terms of theta constants with the named modulus \cite{A2, TM}) and $l_\tau$ is a  linear-fractional map:
\be
\label{norm}
l_\tau(w)=-\frac1{\sqrt{k}}\frac{\sqrt{k_0}w-1}{\sqrt{k_0}w+1},
\qquad\qquad \frac1{\sqrt{k}}=\frac{1+\sqrt{k_0}}{1-\sqrt{k_0}}
\ee 
where the square root of the Jacobi's modulus $\sqrt{k}(\tau)=\theta_2/\theta_3$ is holomorphic in the upper half-plane 
and $\sqrt{k_0}:=\sqrt{k}(\tau_0)$.  

To prove (\ref{ident}) we note that $cd(\dots)$ functions on both sides of the equality are even, degree 2 with the common period lattice 
$Span_{\mathbb{Z}}(2,\tau)$. Hence they differ by a linear fractional map  which may be reconstructed if we evaluate it at three points.
Taking $u=\frac12, ~\frac12+\frac{\tau}2, ~0$ we see that $l_\tau$ sends  three points $(k_0)^{-1/2}, -(k_0)^{-1/2}, 1$
to respectively $0, ~\infty, ~1$  where from the expression (\ref{norm}) for the linear fractional map follows.

Let us now make a transformation of Chebyshev-Blashke product:
$$
y_{n,\tau}(u):=\sqrt{k(n\tau)}cd(2K(n\tau)nu|n\tau)=
\sqrt{k(n\tau)}l_{n\tau}\circ cd(iK'(\tau_1)nu|\tau_1)=
$$ 
(keep in mind the notations: $\tau_1:=-4/(n\tau)$; $l^0_{\tau}:=\sqrt{k(\tau)}l_{\tau}$ and the identity \cite{A2} $cd(u|\tau)=sn(u+K(\tau)|\tau)$),
$$
=l^0_{n\tau}\circ cd(4K(\tau_1)u/\tau|\tau_1)=
l^0_{n\tau}\circ sn(K(\tau_1)(4u/\tau+1)|\tau_1)=
l^0_{n\tau}\circ x_{\tau_1}(4u/\tau+1).
$$
For $n=1$ we in particular have 
$$
y_{1,\tau}(u):=l^0_{\tau}\circ x_{\tau_0}(4u/\tau+1)
$$
with $\tau_0=n\tau_1=-4/\tau$ defined in the above passage. Introducing $1+4u/\tau$ as a new independent variable in the plane
we see that Chebyshev-Blaschke product is just a renormalization of Zolotar\"ev fraction:
$$
Y_n(\cdot|\tau):=y_{n,\tau}\circ y^{-1}_{1,\tau}=l^0_{n\tau}\circ x_{\tau_1}\circ x_{n\tau_1}^{-1}\circ(l^0_{\tau})^{-1}=
l^0_{n\tau}\circ Z_n(\cdot|\tau_1)\circ(l^0_{\tau})^{-1}.
$$

Now it would be usefull to obtain the numerous properties \cite{Ng} of  $Y_n$ just transferring them from the properties of 
Zolotar\"ev fractions.

\section{Accounting of the (Sub)lattices }
We start with enumeration of given index $n$ sublattices  of a rank 2 lattice. Fix any 'standard' basis in $L$, a basis in the 
sublattice $M$ is related to the former by an integer $2\times2$ matrix $Q$, $|det Q|=|L:M|=n$. The change of basis in the sublattice results in multiplying $Q$ on the left by an integer matrix with determinant $\pm1$.
Thus all matrices $Q$ are collected into classes of equivalence representing the same sublattice $M$. 
It is easy to check that each class contains a unique matrix of the kind
$$
Q=
\begin{array}{||cc||}
p&l\\
0&q
\end{array},
\qquad p,q>0; \quad q>l\ge0; \quad pq=n, 
$$
known as a Hermitian normal form of integer matrix \cite{Newman}. 

\begin {rmk}
The representation of a sublattice $M$ by a matrix as above depends on the choice of the 'standard' basis in the ambient lattice
$L$. In other words, there is a right action of invertable integer matrices on Hermitian normal forms. 
\end {rmk}

Each positive integer divisor $q$ of index $n$ gives us $q$ sublattices which differ by the value of 
off-diagonal element of matrix $Q$, therefore we arrive to 
\begin{lmm}
Rank two lattice has $\sigma_1(n)$ sublattices of a given index $n$,
$\sigma_1(n)$ is the sum of all positive integer divisors of $n$.
\end{lmm} 

\begin{xmpl} 
$\sigma_1(11)=1+11=12$; $\sigma_1(12)=1+2+3+4+6+12=28$.
\end{xmpl}

The above Theorem \ref{ZfrEqLat} reduces counting classes of degree $n$ twisted Zolotar\"ev fractions to the enumeration of index $n$  sublattices of a given lattice modulo dilations and rotations. 

\subsection{General lattice} When $j(B)\neq 0,1$ the lattice $L$ has no symmetries:  
$L=\mu L$ means that $\mu=\mp1$. Each index $n$ sublattice $M$ of $L$ gives a class of degree $n$ Zolotar\"ev fractions,
whose number is therefore $\sigma_1(n)$. In case $n=4$ the answer is one less as we throw away the exceptional sublattice 
$M=2L$.

\subsection{Square lattice} When $j(D)=1$ the corresponding lattice $L$ is equivalent to Gaussian integers $\mathbb{Z}[i]$. It survives after multiplication by $i$, therefore the sublattice $M\subset L$ will be in the same class as $iM$ and will be counted twice unless $M=iM$. Each square (sub)lattice $M$ is uniquely determined by its nonzero element $e=p+iq$ of minimal length lying in the sector $\{0\le Arg(u)<\pi/2 \}$. This sublattice has index 
\be
(|e|^2):=n=p^2+q^2, 
\qquad p>0; ~q\ge0.
\label{square}
\ee
Given $n$, the number of ordered integer solutions $p,q$ of the above equation (\ref{square}) is known as $S_2(n)$,
for instance $S_2(25)=3$ since $25=5^2+0^2=4^2+3^2=3^2+4^2$. In this notation,
the number of classes of degree $n$ twisted Zolotar\"ev fractions with $j=1$ is $\frac12(\sigma_1(n)+S_2(n))$.

\subsection{Hexagonal lattice} When $j(D)=0$ the corresponding lattice $L$ is equivalent to $\mathbb{Z}[\epsilon]$, 
$\epsilon=\exp(i\pi/3)$. It survives after the rotation by $\pi/3$. The sublattice $M$ will be in the same class as $\epsilon M$
and $\epsilon^2 M$ and will be counted thrice unless the sublattice $M$ is hexagonal itself. Each hexagonal (sub)lattice $M$ is uniquely determined by its nonzero element $e=p+\epsilon q$ of minimal length lying in the sector $\{0\le Arg(u)<\pi/3 \}$. This sublattice has index 
\be
(|e|^2):=n=p^2+pq+q^2, 
\qquad p>0; ~q\ge0.
\label{hex}
\ee
Given $n$, the number of ordered integer solutions $p,q$ of equation (\ref{hex}) we designate as $h(n)$,
for instance $h(13)=2$ since $13=1^2+1\cdot3+3^2=3^2+3\cdot1+1$. In this notation,
the number of classes of degree $n$ Zolotar\"ev fractions with $j=0$ is $\frac13(\sigma_1(n)+2h(n))$. 

The results of this section may be summarized in a table:

\begin{table}[h!]
\begin{tabular}{c|c|c|c} 
&$j=0$&$j=1$&$j\neq 0,1$\\
\hline
$n=4$& $2$& $3$& $6$\\
$4\neq n>2$& $\frac13(\sigma_1(n)+2h(n))$ & $\frac12(\sigma_1(n)+S_2(n))$& $\sigma_1(n)$\\
\end{tabular}
\caption{Number of classes of twisted Zolotar\"ev fractions, given $j$ invariant and degree $n$}
\end{table}

\parbox{9cm}
{\it
119991 Russia, Moscow GSP-1, ul. Gubkina 8,\\
Institute for Numerical Mathematics,\\
Russian Academy of Sciences\\[3mm]
{\tt ab.bogatyrev@gmail.com}}
\end{document}